\documentstyle[11pt]{article}
\title{{\bf Comment on a Paper by Yucai Su On Jacobian Conjecture}}
\author{{\Large T.T. Moh}\thanks{ Math Department, Purdue University, West Lafayette, Indiana 47907-1395. tel: (765)-494-1930, e-mail ttm@math.purdue.edu}}
\begin{document}
\maketitle

\begin{abstract}
The said paper [2] entitled "Proof Of Two Dimensional Jacobian Conjecture" is with gaps.
\end{abstract}

\noindent{\bf Comments}

J. Lipman sent us the URL of [2]. Harley Flanders of Michigan mentioned that: 
"{\it I would be very much interested in your opinions of Su's paper.}" 

We have browsed Su's paper. There are at least $2$ errors. (1): The author relied on the computations over 
a non-existing entity, {\it the associative algebra {\bf F}[x,y][[$y^{-1}$]]} (cf pg 5 of [2]).  It is easy 
to see that if such a thing exists, we must have
$$
1+yy^{-1}+\cdots+y^iy^{-i}+\cdots=1+1+\cdots+1+\cdots
$$
This is impossible. What the author really means is the ring $K[x]((y^{-1}))$, the ring of meromorphic functions
over $K[x]$. The author should make the substitution throughout his paper.

(2): The author first use an (tame) automorphism $x\mapsto x+y^k,y\mapsto y$ to satisfy the condition that
$f_0(x),g_0(x)\in F^{*}$ (on pg $7$) (note that the degrees are increased), then use the same kind 
automorphism to reduce the degrees (on pg $15$). This is a cyclic argument. The degrees may not be reduced.  
 
\vskip 0.2 in

\noindent{\bf Suggestions to Yucai Su}

We think that the author is making a noble attempt to solve a hard problem. The author should make the substitution 
mentioned in (1) above and consult article [1] (in there, $\tau=F(x,y)^{-1/m}$ in Su's notation). To avoid a cyclic 
argument, one should use a general linear transformation instead of an (tame) automorphism $x\mapsto x+y^k,y\mapsto y$ 
at the beginning. For the last thirty years, the only case one can not handle is the case of "more than one point at 
$\infty$", i.e., in Su's notations, $p=1, j=m+n-2$ (cf pg 17, Subcase), all other cases are well-known. The author 
should give an convincing argument for this case (in fact, for this case only).

\vspace*{1cm}
\noindent{\bf Reference} 

[1] {\bf Moh, T.T.}:  {\it On the Jacobian Conjecture} Crelle 1983.

[2] {\bf Su, Yucai}  {\it Proof Of Two Dimensional Jacobian Conjecture} Arxiv, Dec 18, 2005,
\end{document}